\documentclass[11 pt]{amsart}
\usepackage{graphicx}
\usepackage[hidelinks]{hyperref}
\usepackage{color}
\usepackage{amssymb,amsmath,amsfonts,latexsym}
\usepackage{epstopdf}
\usepackage{nicefrac}
\usepackage{amsthm}
\usepackage{hyperref}
\hypersetup{
     colorlinks   = true,
     citecolor    = blue
}

\usepackage{tikz}

\usepackage{tikz,fullpage}
\usetikzlibrary{patterns}
\usepackage{enumerate}
\usepackage{ulem}
\usepackage{float, lscape}
\numberwithin{equation}{section}
\usepackage[margin=1.27in]{geometry}
\usetikzlibrary{arrows,%
                petri,%
                topaths}%
\usepackage{tkz-berge}
\usepackage[position=top]{subfig}
\usetikzlibrary{%
  matrix,%
  calc,%
  arrows%
}
\newtheorem{thm}{THEOREM}[section]

\begin{document}
\title{Equivariant Cartan Homotopy Formulae for the crossed product of $DG$ Algebra}
\author{Safdar Quddus}

\date{\today}
 
\let\thefootnote\relax

\begin{abstract}
We establish the equivariant Cartan Homotopy Formula for the crossed product of $DG$-algebra obtained by a finite group action.
\end{abstract}
\maketitle
\section{Introduction}
The paracyclic modules were used in \cite{GJ1} and \cite{BGJ} to understand the cyclic homology of the crossed product algebras. It is an analogue of the Eilenberg-Zilber theorem for bi-paracyclic modules. Let $A$ be a unital $DG$ algebra over a commutative ring $k$ and let $G$ be a finite discrete group which acts on $A$ by automorphisms. The result is that the cyclic homology of the crossed product algebra $A \rtimes G$ has the following decomposition.
\begin{thm}[\cite{GJ1} \cite{FT}]
If $G$ is finite and $|G|$ is invertible in $k$, then there is a natural isomorphism of cyclic homology and 
\begin{center}
$$HC_\bullet (A\rtimes G;W) = HC_\bullet(H_0(G, A^\natural_G);W)$$,
\end{center} 
where $H_0(G, A^\natural_G)$ is the cyclic module
\begin{center}
$$ H_0(G, A^\natural_G)(n) = H_0 (G,k[G]\otimes A^{(n+1)}).$$
\end{center} 
\end{thm}
Where $W$ is a finite dimensional graded module over the polynomial ring $k[u]$, where $deg(u) = -2$; the above result when considered for different coefficients $W$ yield several theories, some are illustrated below:\\

$1)$  $W= k[u]$ gives negative cyclic homology $HC^{-}(A)$;\\
$2)$  $W = k [u, u^{-1}]$ gives periodic cyclic homology $HP_\bullet (A)$;\\
$3)$ $W=k[u, u^{-1}]/ uk[u]$ gives cyclic homology $HC_\bullet(A)$;\\
$4)$  $W= k[u]/uk[u]$ gives the Hochschild homology $HH_\bullet(A)$.\par
The above decomposition has been studied and used to understand the homological properties of crossed product algebras \cite{CGGV} \cite{P1} \cite{P2} \cite{P3} \cite{Q1} \cite{Q2}\cite{Q3} \cite{Q4} \cite{ZH}. The prima fact is that the (co)homology modules of the crossed product algebra decomposes into twisted (co)homology modules relative to the conjugacy classes of $G$. Hence the studying the crossed product algebra is simplified.
\par
\subsection{$DG$ Algebra}
A unital DG algebra $(A, d)$ is a differentially graded unital algebra $A$, with $k$-bilinear maps
\begin{center} 
$A_n \times A_m \to A_{n+m}$,
 sending $(a,b) \mapsto ab$ such that
 \end{center}
 \begin{center}
 $d_{n+m}(ab)=d_n(a)b + (-1)^n a d_m(b)$
 \end{center}
and such that $ \oplus A_n$ becomes an associative and unital $k$-algebra. Through out this article we demand that the finite group action on $(A, d)$ preserves the grading and commutes with the differential structure; $i.e.$ for $g \in G$, $d g = g d$. The (co)homology theories on $DG$ algebras have been studied extensively \cite{T} \cite{K} \cite{GJ2}. The above algebra can also be considered as a special $A_\infty$-algebra $(A,m_i)$ with $m_i=0$ for $i>2$. Here $m_2$~(2-co-chain for $A$) defines product on $A$ and $m_1$ is the mixed differential.

\subsection{Cartan Homotopy Formulae}
The Cartan homotopy formulae for algebra were first observed by Rinehart \cite{R} in the case where $D$ is a derivation on commutative algebra, and later in full generality by Getzler \cite{G} for $A_\infty$-algebras. The formulae is stated below:

Given $D \in C^k(A,A)$, we consider its Lie derivative $L_D$ and the Hochschild and cyclic chain maps $b$ and $B$. 
Let $\iota_D$ denote the contraction associated to $D$ and $S_D$ be the corresponding suspension.

\begin{thm}[Cartan Homotopy Formulae]
$$[b, L_D]+L_{\delta D}=0, [B,L_D]=0 \text{ and } [L_D, L_E]=L_{[D,E]}$$ 
$$[b+B, \iota_D + S_D] = L_D + \iota_{\delta D} + S_{\delta D}$$
\end{thm}
where $[,]$ is the graded Gerstenhaber bracket for $C^\bullet(A,A)$.\\
The above formulae was proved in generality for the $A_\infty$-algebras by Getzler \cite{G}. The case $k=1$ in his paper corresponds to the DG-algebra $(A,d)$. Given a crossed product algebra we can ask if the Cartan homotopy formulae hold for each components in the decomposition described above. As to my knowledge, the question is unanswered even for an associative algebra. We shall answer this question for $DG$-algebras when the group action preserves the gradation.\\

\section{Statement}
For $G$ a finite discrete group acting on a $DG$-algebra $(A,d)$ over a ring $k$ such that $|G|$ is invertible in $k$. For $g \in G$ define a $g$-twisted $A$-left-module structure on $A$ by the following formula:

$$a\bullet (u_{g} m):= u_{g} g^{-1}(a)\cdot m,$$ 
where $a,m \in A$. We tag by $u_g$ the twisted element of $A$ with this left bi-module structure and the left twisted bi-module is denoted by $A_g$~(refer \cite{Q1}, pp 332).

Define 
\[
C_{0}(A)_g:=A_{g} \quad \text { and } \quad C_{n}(A)_{g}:=A_{g} \otimes A^{\otimes n}.
\]

Let $b^g$, $B^g$ denote the chain differentials of the complex $C_\bullet(A)_g$ and $L_D^g$ be the twisted Lie derivative associated to $D \in C^k(A,A)$. We shall produce the explicit expression of these maps in this literature.
\begin{thm}
The Lie derivative $L_D^g$, chain maps $b^g$ and $B^g$ satisfy the following:
\[
[b^g , L_D^g]+ L_{\delta D}^g=0 \quad [B^g , L_D^g]=0 \quad [L_D^g, L^g_E] = L_{[D,E]}^g \quad \text { and } [b^g+B^g, \iota_D^g + S_D^g] = L_D^g + \iota^g_{\delta D} + S_{\delta D}^g.
\]

\end{thm}
The $[,]$ above is the Gerstenhaber bracket for $C^\bullet(A,A)$ and $\delta$ is the  Hochschild co-chain map.
\section{Paracyclic decomposition for cross product $(A,d) \rtimes G$}
\subsection{Hochschild Chain Complex for $(A,d)$} 
Define the differentials $d : C_\bullet(A) \to C_\bullet(A)$, $b : C_\bullet(A) \to C_{\bullet}(A)[-1]$ and $B : C_\bullet(A) \to C_{\bullet}(A)[1]$ as follows.\\
\begin{center}
$$d(a_0\otimes \cdots \otimes a_n)=\sum_{i=1}^n (-1)^{\sum_{k<i}(|a_k|+1)+1}(a_0 \otimes \cdots \otimes da_i \otimes \cdots \otimes a_n)$$
\end{center}
\begin{center}
$$b(a_0 \otimes \cdots \otimes a_n) = \sum_{k=0}^{n-1} (-1)^{\sum_{i=0}^{k}(|a_i|+1)+1}(a_0 \otimes \cdots \otimes a_k a_{k+1} \otimes \cdots \otimes a_n)$$ \\ $$+(-1)^{{|a_n|} +(|a_n|+1){\sum_{i=0}^{n-1}(|a_i|+1)+1}}(a_na_0 \otimes \cdots  \otimes a_{n-1})$$
\end{center} and 

\begin{center}
$$B(a_0 \otimes \cdots \otimes a_n)=\sum_{k=0}^n(-1)^{\sum_{i\leq k}{(|a_i|+1){\sum_{i \geq k}(|a_i|+1)}}}(e \otimes a_{k+1} \otimes \cdots \otimes a_n \otimes a_0 \otimes \cdots  \otimes a_{k})$$
\end{center}

Where $|a|$ denotes the degree of $a$, also the above formula satisfies $B(e, a_1, . . . , a_k) = 0$. These sign conventions ensures that the elements $a_i$ for $i > 0$ occur with an implicit suspension reducing the degree to $|a_i| -1$. The homology of $(C(_\bullet(A, A), b+d)$ is the Hochschild homology of $(A,d)$ with coefficients in $A$ endowed with the bi-module structure. Let $W$ be a graded module over the polynomial ring $k[u]$, where $deg(u) = -2$; such that $W$ has finite homological dimension then the module $H_\bullet(C_\bullet \boxtimes W, b+d+uB)$, where $(C_\bullet \boxtimes W)= C_\bullet[[u]] \otimes_{k[u]} W$, is the cyclic homology of the mixed complex $(A, b+d, B)$ with coefficients in $W$. Some examples of cyclic (co)homology for various $W$ are listed the statement of Theorem 0.1.

\section{Proof of the theorem}
The paracyclic decomposition of the cyclic homology of $(A,d) \rtimes G$ exists and the following is the decomposition: 
$$HC_\bullet(A\rtimes G)= \bigoplus_{[g]} HC_\bullet(H_0(G^g, A^\natural_g))=\bigoplus_{[g]} HC_\bullet(A_g)^G$$ with the maps; 
$b^g : C_\bullet(A)_g \to C_{\bullet}(A)_g[-1]$ and $B^g : C_\bullet(A)_g \to C_{\bullet}(A)_g[1]$ defined as:

$$b^g(u_ga_0 \otimes \cdots \otimes a_n) = \sum_{k=0}^{n-1} (-1)^{\sum_{i=0}^{k}(|a_i|+1)+1}(u_g a_0 \otimes \cdots \otimes a_k a_{k+1} \otimes \cdots \otimes a_n)$$ \\ $$+(-1)^{{|a_n|} +(|a_n|+1){\sum_{i=0}^{n-1}(|a_i|+1)+1}}(u_g g^{-1}(a_n)a_0 \otimes \cdots  \otimes a_{n-1})$$
and $$B^g(u_g a_0 \otimes \cdots \otimes a_n)=\sum_{k=0}^n(-1)^{\sum_{i\leq k}{(|a_i|+1){\sum_{i \geq k}(|a_i|+1)}}}(u_g e \otimes g^{-1}(a_{k+1}) \otimes \cdots \otimes g^{-1}(a_n) \otimes a_0 \otimes \cdots  \otimes a_{k})$$
To explicitly derive the above maps we consider the definition of these in terms of the elementary maps $s_i$, $t$ and $d_i$ \cite{L} and decipher the twisted elementary maps \cite{BG} and finally formulate $b^g$ and $B^g$.
One important point here to be noted is that since the action of $G$ on $(A,d)$ preserves grading, the crucial sign term $\epsilon_k = {\sum_{i\leq k}{(|a_i|+1){\sum_{i \geq k}(|a_i|+1)}}}$ remains unchanged~(in fact under this assumption the signs for $A_\infty$ -algebras are also invariant under the unified algebraic structure for chains and co-chains of $A_\infty$-algebras ) and is a necessary condition to yield the Cartan homotopy formulae for twisted components. \par
For $D \in C^k(A,A)$, the maps $\iota_D$ and $S_D^g$ can be computed using the signed cyclic permutation map $t$ \cite{BGJ}. They are as follows:

$$\iota_D(u_g a_0 \otimes \cdots \otimes a_n)= (-1)^{|D||a_0|} a_0 D(a_1, \dots , a_n);$$

$$L_D^g(u_g a_0 \otimes \cdots \otimes a_n) = \sum_{k=1}^{n-d} (-1)^{\nu_k(D, n)} u_g a_0 \otimes \dots \otimes D(a_{k+1}, \dots , a_{k+d}) \otimes \dots a_n$$\\
$$+\sum_{k=n+1-d}^n (-1)^{\eta_k(D,n)} u_g D(g^{-1}(a_{k+1}), \dots , g^{-1}(a_n), a_0, \dots)\otimes \dots \otimes a_k.$$

where the terms inside $D$ in the second summand must contain $a_0$~(without $g$ action) and be cyclically permuted.

$$S_D(u_g a_0 \otimes \cdots \otimes a_n)=\sum_{j\geq 0; k\geq j+d} (-1)^{\epsilon_{jk}(D,n)} u_g e \otimes g^{-1}(a_{k+1}) \otimes \dots g^{-1}(a_{n})\otimes a_0 \otimes \dots \otimes D(a_{j+1}, \dots , a_{j+d})\otimes \dots \otimes a_k.$$

where, $|D|=\text{ (degree of the linear map D) } + d$, $D \in C^d(A,A)$ is being considered as a linear map $D : A^{\otimes d} \to A$; and the sign coefficients are $\nu_k(D,n) = (|D|+1)(|a_0| + \sum_{i=1}^{k}(|a_i|+1))$, $\eta_k(D,n) = |D| + \sum_{i\leq k}(|a_i| +1) \sum_{i \geq k}(|a_i|+1)$ and  $$\epsilon_{jk}(D,n)=(|D|+1)(\sum_{i=k+1}^n (|a_i|+1) +|a_0|+\sum_{i=1}^j(|a_i|+1)).$$

We briefly describe the Gerstenhaber algebra structure on the co-homology $H^\bullet(A,A)$, the cup product is defined as below, for $D \in C^d(A,A)$ and $E \in C^e(A,A)$. 

$$(D \smile E)(a_1, \dots , a_{d+e}) = (-1)^{|E| \sum_{i \leq e}(|a_i| +1)} D(a_1, \dots , a_e)E(a_{d+1}, \dots, a_{d+e});$$
and the product $\circ$ is defined as:
$$(D\circ E)(a_1, \dots , a_{d+e})= \sum_{j\geq 0}(-1)^{(|E|+1)\sum_{i=1}^j(|a_i|+1)}D(a_1,\dots, a_j, E(a_{j+1}, \dots, a_{j+e}), \dots).$$
The Lie bracket is hence $[D,E] = D \circ E - (-1)^{|D|+1)(|E|+1)|} E \circ D$. The co-chain map $\delta$ on $C^\bullet(A,A)$ is 

$$(\delta D)(a_1, \dots, a_{d+1}) = (-1)^{|a_1||D| + |D| +1} a_1 D(a_2, \dots , a_{d+1})$$
$$+\sum_{j=1}^d (-1)^{|D|+1+\sum_{i=1}^j(|a_i|+1)} D(a_1, \dots , a_{j}a_{j+1}, \dots, a_{d+1}) + (-1)^{|D|\sum_{i=1}^d(|a_i| +1)}D(a_1, \dots, a_d)a_{d+1}$$
The operator $\delta D$ can also be described as $\delta D = [m_2,D]$.
\begin{proof}[Proof of Theorem 1.1]
\begin{equation}
[b^g, L_D^g]+L^g_{\delta D}=b^g L_D^g - L_D^g b^g+L^g_{\delta D}. \label{eq:1}
\end{equation}
We evaluate the above expression on the element $(u_g a_0 \otimes \cdots \otimes a_n)$. We define the bar complex map on the equivariant Hochschild chain complex by $b^{g '}$, it is as follows:

$$b^{g '}(u_g a_0 \otimes \cdots \otimes a_n)=\sum_{k=0}^{n-1} (-1)^{\sum_{i=0}^{k}(|a_i|+1)+1}(u_g a_0 \otimes \cdots \otimes a_k a_{k+1} \otimes \cdots \otimes a_n)$$

Similarly we define the operator $L_D^{g '}$ as follows:
$$L_D^{g '}(u_g a_0 \otimes \cdots \otimes a_n)=\sum_{k=1}^{n-d} (-1)^{\nu_k(D, n)} u_g a_0 \otimes \dots \otimes D(a_{k+1}, \dots , a_{k+d}) \otimes \dots a_n.$$

The untwisted expression in  \eqref{eq:1} is $b^{g '}L_D^{g'} - L_D^{g '}b^{g '} + L_{\delta D}^{g '}$. We collect the coefficients and the terms cancel each other as signs mismatch, for example the expression $u_g a_0 \otimes \dots \otimes a_w D(a_{w+1}, \cdots, a_{w+d})\otimes \dots a_n$ cancel only if
$$\nu_{w}(D, n)+ \sum_{i=0}^w (|a_i|+1) +1+ \nu_{w-1}(\delta D, n)+ |a_w||D| +|D|+1 \equiv 1(\mod 2).$$
Which is true as $|\delta D| = |D|+1$. The relation \eqref{eq:1} is $0$ in the untwisted case hence the untwisted terms cancel each other once the signs mismatch. On the other hand the cancellation of the twisted terms in \eqref{eq:1} is non-trivial. To see this we firstly observe that
$$|D(a_{i+1}, \dots ,a_{i+d})| = |D|\sum_{j=i+1}^{i+d}(|a_j|+1) \text{ and } |g^{-1}(a)| = |a|$$
Using the above relations we can see that the parity of signs remain the same and terms cancel out as it did in the untwisted. For example, \\
$$\eta_k(\delta D, n) + |D| \sum _{i=1}^{d}(|b_i|+1) + \eta_k(D,n) + |D(b_1, \dots , b_d)|+1\equiv 1(\mod 2), $$
where $(b_1, \dots , b_d) = (g^{-1}(a_{k+1}), \dots , g^{-1}(a_n), a_0, \dots , a_{d-n+k})$ and hence the terms of kind $u_g D(g^{-1}(a_{k+1}),\dots g^{-1}(a_n), a_0, \dots )a_{d-n+k+1}\otimes \cdots \otimes a_k$ cancel each other. Similarly, it is easy to check that all other types of twisted terms cancel each other and the appropriate parity of signs are ensured by the grade preserving group action.

\begin{equation} \label{eq:2}
[L_D^g, L_E^g]=L_{[D,E]}^g
\end{equation}
The proof of the above relation is straight forward, the bracket $[,]$ is the Gerstenhaber Lie algebra commutator as described above.

\begin{equation} \label{eq:3}
[B^g, L_D^g] = B^g L_D^g - L_D^g B^g.
\end{equation}
$$B^g L_D^g(u_g a_0 \otimes \cdots \otimes a_n)=B^g \Big\{ \sum_{l=1}^{n-d} (-1)^{\nu_l(D, n)} u_g a_0 \otimes \dots \otimes D(a_{l+1}, \dots , a_{l+d}) \otimes \dots a_n$$\\
$$+\sum_{l=n+1-d}^n (-1)^{\eta_l(D,n)} u_g D(g^{-1}(a_{l+1}), \dots , g^{-1}(a_n), a_0, \dots)\otimes \dots \otimes a_l \Big\}$$
$$=\sum_{l=1}^{n-d} (-1)^{\nu_l(D, n)} \sum_{k=0}^{n-d+1}(-1)^{\sum_{i\leq k}{(|a'_i|+1){\sum_{i \geq k}(|a'_i|+1)}}} u_g e \otimes g^{-1}(a'_{k+1}) \otimes \cdots \otimes g^{-1}(a'_n) \otimes a'_0 \cdots \otimes a_k$$
Such that any of $a'_i$ could be $D(\dots)$.
$$+\sum_{l=n+1-d}^n (-1)^{\eta_l(D,n)} \sum_{k=0}^{n-d+1}(-1)^{\sum_{i\leq k}{(|a''_k|+1){\sum_{i \geq k}(|a''_k|+1)}}} u_g e \otimes g^{-1}(a_{k+1})\otimes \cdots D(g^{-1}(a_{l+1}), \dots , g^{-1}(a_n), a_0, \dots)\otimes $$
$$\dots \otimes a_k $$
such that for each $l$ one of the $a''_k$ in the sign expression is $D(g^{-1}(a_{l+1}), \dots , g^{-1}(a_n), a_0, \dots)$ and the rest are the remaining $a_i $'s that do not appear in the expression $D(g^{-1}(a_{l+1}), \dots , g^{-1}(a_n), a_0, \dots)$.\\

$$L_D^g B^g(u_g a_0 \otimes \cdots \otimes a_n)=L_D^g \Big\{ \sum_{k=0}^n(-1)^{\sum_{i\leq k}{(|a_i|+1){\sum_{i \geq k}(|a_i|+1)}}}(u_g e \otimes g^{-1}(a_{k+1}) \otimes \cdots \otimes g^{-1}(a_n) \otimes a_0 \otimes \cdots  \otimes a_{k}) \Big\}$$
$$=L_D^{g '}\Big\{ \sum_{k=0}^n(-1)^{\sum_{i\leq k}{(|a_i|+1){\sum_{i \geq k}(|a_i|+1)}}}(u_g e \otimes g^{-1}(a_{k+1}) \otimes \cdots \otimes g^{-1}(a_n) \otimes a_0 \otimes \cdots  \otimes a_{k}) \Big\}+0$$
$$=\sum_{k=0}^n\sum_{l = 1}^{n+1-d}(-1)^{\nu_l (D, n+1)+\sum_{i\leq k}{(|a_i|+1){\sum_{i \geq k}(|a_i|+1)}}}u_g e \otimes \cdots \otimes D(g^{-1}(a_{l+1}), \dots, g^{-1}(a_n), a_0,...)\otimes \cdots a_{k}.$$
In the above expressions we have repeatedly used the fact that $|g^{-1}(a)|=|a|$. We also observe that $|D(g^{-1}(a_{l+1}), \dots, g^{-1}(a_n), a_0,...)|=|D|(\sum_{j=0}^{d-n+l}(|a_j|+1)+\sum_{j=l+1}^{n}(|a_j|+1))$. Hence the sign mismatch for cancellation to yield that $L_D^g$ and $B^g$ commute.

Finally we want to show that:
\begin{equation} \label{eq:4}
[b^g+B^g, \iota_D^g + S_D^g] = L_D^g + \iota^g_{\delta D} + S_{\delta D}^g
\end{equation}
The expression above can be written as below:
$$[b^g, \iota_D^g] + [B^g , \iota_D^g] + [b^g , S_D^g] + [B^g , S_D^g] = P + Q + R + 0 $$
Observe that $[b^g, \iota^g_D]=\iota^g_{\delta D}$ because the twisted terms $(b^g-b^{g '})\iota^g_D$ and  $\iota^g_D (b^g-b^{g '})$ cancel each other. Hence we are left to show that:
$$Q+R = L_D^g + S_{\delta D}^g .$$
The above relation can be seen be comparing the parity of sign indices for the terms, for example the sign for the term $u_g D(g^{-1}(a_k), \dots , g^{-1}(a_n), \dots, ) \otimes  \cdots \otimes a_k)$ for the RHS is $\eta_k(D,n)$ while in the LHS it has sign coefficient as $|D||e| + \sum_{i \leq k}(|a_i|+1)\sum_{i \geq k} (|a_i|+1)$; but indeed 
$$ \eta_k(D,n) \equiv |D||e| + \sum_{i \leq k}(|a_i|+1)\sum_{i \geq k} (|a_i|+1) \mod 2$$
since $|e| =1$. Hence the terms of the given form cancel out. Similar relations involving $\eta_k(D,n)$, $\nu_k(D,n)$, $\epsilon_{jk}(D,n)$and $\epsilon_{jk}(\delta D,n)$ yield the desired result.
\end{proof}

\section{Acknowledgement}

This result grew out of the discussion and collaboration with Xiang Tang and Sayan Chakraborty on equivariant Gauss-Manin connections, and the discussion has contributed to the outcome of the paper.

\vspace{2mm}

Email: safdarquddus@iisc.ac.in.


\begin{thebibliography}{}


\bibitem[BG]{BG} Block, J.; Getzler, E. Equivariant cyclic homology and equivariant differential forms. Annales scientifiques de lEcole Normale Superieure, Serie 4, Volume 27 (1994) no. 4, pp. 493-527.


\bibitem[BGJ]{BGJ} Block, J.; Getzler, E.; Jones, J. D. S. The cyclic homology of crossed product algebras. II. Topological algebras. J. Reine Angew. Math. 466 (1995), 19-25.

\bibitem[CGGV]{CGGV}  Carboni, G.; Guccione, J. A.; Guccione, J. J.; Valqui, C.; Cyclic homology of Brzeziński's crossed products and of braided Hopf crossed products. Adv. Math. 231 (2012), no. 6, 3502-3568

\bibitem[FT]{FT} Feigin, B.L. and Tsygan B., Additive K-theory, In “K-Theory, Arithmetic and
Geometry”, Lecture Notes in Mathematics 1289 (Yu. I. Manin, ed.), 1986, pp. 67-209.

\bibitem[G]{G} Getzler, E.; Cartan homotopy formulas and the Gauss-Manin connection in cyclic homology. Quantum deformations of algebras and their representations (Ramat-Gan, 1991/1992; Rehovot, 1991/1992), 65–78, Israel Math. Conf. Proc., 7, Bar-Ilan Univ., Ramat Gan, 1993.

\bibitem[GJ1]{GJ1} Getzler, E.; Jones, J. D. S. The cyclic homology of crossed product algebras. J. Reine Angew. Math. 445 (1993), 161-174

\bibitem[GJ2]{GJ2} Getzler, E., and J. D. S. Jones. "$A_\infty$-algebras and the cyclic bar complex." Illinois Journal of Mathematics 34, no. 2 (1990): 256-283

\bibitem[L]{L} J. Loday, Cyclic Homology, second edition, Springer, ISBN 3540630740, 1998.

\bibitem[K]{K} Khalkhali, M.; On Cartan homotopy formulas in cyclic homology, Manuscripta Math. 94 (1) (1997) 111-132.
 
\bibitem [P1]{P1}  Ponge, R.; The cyclic homology of crossed-product algebras, I. C. R. Math. Acad. Sci. Paris 355 (2017), no. 6, 618-622.

\bibitem [P2]{P2}  Ponge, R.; The cyclic homology of crossed-product algebras, II. C. R. Math. Acad. Sci. Paris 355 (2017), no. 6, 623-627.

\bibitem [P3]{P3}  Ponge, R.; Cyclic homology and group actions. J. Geom. Phys. 123 (2018), 30-52.

\bibitem[R]{R} Rinehart, G. S.; Differential forms on general commutative algebras, Trans. Amer. Math. Soc., 108 (1963), 195-222.

\bibitem [Q1]{Q1} Quddus, S.; Hochschild and cyclic homology of the crossed product of algebraic irrational rotational algebra by finite subgroups of $SL(2, \mathbb Z)$. J. Algebra 447 (2016), 322-366. 

\bibitem [Q2]{Q2} Quddus, S.; Cyclic cohomology and Chern-Connes pairing of some crossed product algebras. J. Algebra 481 (2017), 120-157.

\bibitem [Q3]{Q3} Quddus, S.; Cohomology of $\mathcal A_\theta^{alg} \rtimes \mathbb Z_2$ and its Chern-Connes pairing. J. Noncommut. Geom. 11 (2017), no. 3, 827-843.

\bibitem [Q4]{Q4} Quddus, S.; Invariants of the $\mathbb Z_2$ orbifolds of the Podleś' two spheres. J. Noncommut. Geom. 13 (2019), no. 1, 257-267.

\bibitem [T]{T} Tsygan, B.; On the Gauss–Manin connection in cyclic homology, Methods Funct. Anal. Topology, 13 (2007), no. 1, 83-94.

\bibitem[ZH]{ZH}  Zhang, J.; Hu, N.; Cyclic homology of strong cross product algebras. J. Reine Angew. Math. 663 (2012), 177-207

\end{thebibliography}
\end{document}